\documentclass[12pt]{amsart}
\def\N{{\Bbb N}}
\def\Z{{\Bbb Z}}

\newtheorem{Theorem}{Theorem}[section]

\newtheorem{Lemma}[Theorem]{Lemma}

\theoremstyle{definition}

\newtheorem*{Prob}{Problem}
\theoremstyle{remark}

\begin{document}
\sloppy
\title{On a new class of rigid Coxeter groups}
\author{Tetsuya Hosaka} 
\address{Department of Mathematics, Utsunomiya University, 
Utsunomiya, 321-8505, Japan}
\date{January 5, 2005}
\email{hosaka@cc.utsunomiya-u.ac.jp}
\thanks{Partly supported by the Grant-in-Aid for Scientific Research, 
The Ministry of Education, Culture, Sports, Science and Technology, Japan, 
(No.~15740029).}
\maketitle
\begin{abstract}
In this paper, 
we give a new class of rigid Coxeter groups.
Let $(W,S)$ be a Coxeter system. 
Suppose that 
(0) for each $s,t\in S$ such that $m(s,t)$ is even, $m(s,t)\in\{2\}\cup 4\N$, 
(1) for each $s\neq t\in S$ such that $m(s,t)$ is odd, 
$\{s,t\}$ is a maximal spherical subset of $S$, 
(2) there does not exist 
a three-points subset $\{s,t,u\}\subset S$ such that $m(s,t)$ and $m(t,u)$ are odd, 
and (3) for each $s\neq t\in S$ such that $m(s,t)$ is odd, 
the number of maximal spherical subsets of $S$ 
intersecting with $\{s,t\}$ is at most two, 
where $m(s,t)$ is the order of $st$ in the Coxeter group $W$.
Then we show that the Coxeter group $W$ is rigid.
This is an extension of a result of D.~Radcliffe.
\end{abstract}

\section{Introduction and preliminaries}

The purpose of this paper is to give a new class of rigid Coxeter groups.
A {\it Coxeter group} is a group $W$ having a presentation
$$\langle \,S \, | \, (st)^{m(s,t)}=1 \ \text{for}\ s,t \in S \, \rangle,$$ 
where $S$ is a finite set and 
$m:S \times S \rightarrow \N \cup \{\infty\}$ is a function 
satisfying the following conditions:
\begin{enumerate}
\item[(i)] $m(s,t)=m(t,s)$ for any $s,t \in S$,
\item[(ii)] $m(s,s)=1$ for any $s \in S$, and
\item[(iii)] $m(s,t) \ge 2$ for any $s,t \in S$ such that $s\neq t$.
\end{enumerate}
The pair $(W,S)$ is called a {\it Coxeter system}.
For a Coxeter group $W$, a generating set $S'$ of $W$ is called 
a {\it Coxeter generating set for $W$} if $(W,S')$ is a Coxeter system.
Let $(W,S)$ be a Coxeter system.
For a subset $T \subset S$, 
$W_T$ is defined as the subgroup of $W$ generated by $T$, 
and called a {\it parabolic subgroup}.
A subset $T\subset S$ is called a {\it spherical subset of $S$}, 
if the parabolic subgroup $W_T$ is finite.

Let $(W,S)$ and $(W',S')$ be Coxeter systems. 
Two Coxeter systems $(W,S)$ and $(W',S')$ are 
said to be {\it isomorphic}, 
if there exists a bijection 
$\psi:S\rightarrow S'$ such that 
$$m(s,t)=m'(\psi(s),\psi(t))$$ 
for every $s,t \in S$, where 
$m(s,t)$ and $m'(s',t')$ are the orders of $st$ in $W$ 
and $s't'$ in $W'$, respectively.

In general, a Coxeter group does not always determine 
its Coxeter system up to isomorphism.
Indeed some counter-examples are known (cf.\ \cite{Bo}, \cite{BMMN}).
Here there exists the following natural problem.

\begin{Prob}[\cite{BMMN}, \cite{CD}]
When does a Coxeter group determine its Coxeter system up to isomorphism?
\end{Prob}

A Coxeter group $W$ is said to be {\it rigid}, 
if the Coxeter group $W$ determines 
its Coxeter system up to isomorphism 
(i.e., 
for each Coxeter generating sets $S$ and $S'$ for $W$ 
the Coxeter systems $(W,S)$ and $(W,S')$ are isomorphic).

A Coxeter system $(W,S)$ is said to be {\it even}, 
if $m(s,t)$ is even for all $s\neq t$ in $S$.
Also a Coxeter system $(W,S)$ is said to be {\it strong even}, 
if $m(s,t)\in\{2\}\cup 4\N$ for all $s\neq t$ in $S$.

The following theorem was proved by D.~Radcliffe in \cite{R}.

\begin{Theorem}[\cite{R}]\label{Thm:R}
If $(W,S)$ is a strong even Coxeter system, 
then the Coxeter group $W$ is rigid.
\end{Theorem}

In this paper, we prove the following theorem 
which is an extension of Theorem~\ref{Thm:R} and \cite[Theorem~1.2]{H}.

\begin{Theorem}\label{Thm}
Let $(W,S)$ be a Coxeter system.
Suppose that 
\begin{enumerate}
\item[(0)] for each $s,t\in S$ such that $m(s,t)$ is even, 
$m(s,t)\in\{2\}\cup 4\N$, 
\item[(1)] for each $s\neq t\in S$ such that $m(s,t)$ is odd, 
$\{s,t\}$ is a maximal spherical subset of $S$, 
\item[(2)] there does not exist 
a three-points subset $\{s,t,u\}\subset S$ 
such that $m(s,t)$ and $m(t,u)$ are odd, and
\item[(3)] for each $s\neq t\in S$ such that $m(s,t)$ is odd, 
the number of maximal spherical subsets of $S$ 
intersecting with $\{s,t\}$ is at most two.
\end{enumerate}
Then the Coxeter group $W$ is rigid.
\end{Theorem}

\section{Proof of the theorem}

Let $(W,S)$ be a Coxeter system.
Suppose that 
\begin{enumerate}
\item[(0)] for each $s,t\in S$ such that $m(s,t)$ is even, 
$m(s,t)\in\{2\}\cup 4\N$, 
\item[(1)] for each $s\neq t\in S$ such that $m(s,t)$ is odd, 
$\{s,t\}$ is a maximal spherical subset of $S$, 
\item[(2)] there does not exist 
a three-points subset $\{s,t,u\}\subset S$ 
such that $m(s,t)$ and $m(t,u)$ are odd, and
\item[(3)] for each $s\neq t\in S$ such that $m(s,t)$ is odd, 
the number of maximal spherical subsets of $S$ 
intersecting with $\{s,t\}$ is at most two.
\end{enumerate}
Let $(W',S')$ be a Coxeter system.
We suppose that there exists an isomorphism $\phi:W\rightarrow W'$.
To prove Theorem~\ref{Thm}, we show that the Coxeter systems 
$(W,S)$ and $(W',S')$ are isomorphic.

The following lemma is known.

\begin{Lemma}[cf.\ \cite{BMMN}, \cite{R}]\label{lem2-2}
For each maximal spherical subset $T\subset S$, 
there exists a unique maximal spherical subset $T'\subset S'$ 
such that $\phi(W_T)=w'W'_{T'}{w'}^{-1}$ for some $w'\in W'$.
\end{Lemma}

We first prove the following lemma.

\begin{Lemma}\label{lem3-1}
The Coxeter system $(W',S')$ satisfies the following:
\begin{enumerate}
\item[$(0')$] for each $s',t'\in S'$ such that $m'(s',t')$ is even, 
$m'(s',t')\in\{2\}\cup 4\N$, 
\item[$(1')$] for each $s'\neq t'\in S'$ such that $m'(s',t')$ is odd, 
$\{s',t'\}$ is a maximal spherical subset of $S'$, 
\item[$(2')$] there does not exist 
a three-points subset $\{s',t',u'\}\subset S'$ 
such that $m'(s',t')$ and $m'(t',u')$ are odd, and
\item[$(3')$] for each $s'\neq t'\in S'$ such that $m'(s',t')$ is odd, 
the number of maximal spherical subsets of $S'$ 
intersecting with $\{s',t'\}$ is at most two.
\end{enumerate}
\end{Lemma}

\begin{proof}
Let $s'\neq t'\in S'$.
There exists a maximal spherical subset $T'$ of $S'$ 
such that $\{s',t'\}\subset T'$.
By Lemma~\ref{lem2-2}, 
$\phi^{-1}(W'_{T'})\sim W_T$ for some maximal spherical subset $T$ of $S$.
By $(0)$ and $(1)$, either 
\begin{enumerate}
\item[(i)] $(W_T,T)$ is a strong even Coxeter system, or 
\item[(ii)] $|T|=2$ and if $T=\{s,t\}$ then $m(s,t)$ is odd.
\end{enumerate}
Hence $W_T$ is a rigid Coxeter group and 
$(W_T,T)$ and $(W'_{T'},T')$ are isomorphic.
Thus if $m'(s',t')$ is even then $m'(s',t')\in\{2\}\cup 4\N$, 
and if $m'(s',t')$ is odd then 
$\{s',t'\}$ is a maximal spherical subset of $S'$.
Hence $(0')$ and $(1')$ hold.
We can show $(2')$ and $(3')$ 
by the same argument as the proof of \cite[Lemma~3.1]{H}
\end{proof}

Let $\mathcal{A}$ and $\mathcal{A}'$ be 
the sets of all maximal spherical subsets of $S$ and $S'$, respectively.
For each $T\in\mathcal{A}$, 
there exists a unique element 
$T'\in\mathcal{A}'$ such that $\phi(W_T)\sim W'_{T'}$ 
by Lemma~\ref{lem2-2}.

We define 
\begin{align*}
\bar{S}&=\bigcup\{T\in\mathcal{A}\,|\,(W_T,T) \text{\ is strong even}\} \\
\bar{S'}&=\bigcup\{T'\in\mathcal{A'}\,|\,(W'_{T'},T') \text{\ is strong even}\}.
\end{align*}
We note that 
for each $s\in S\setminus \bar{S}$, 
there exists a unique element $t\in S\setminus\{s\}$ 
such that $m(s,t)$ is odd.
Then $m(s,u)=\infty$ for any $u\in S\setminus\{s,t\}$.

Let $W^{{\rm ab}}$ and ${W'}^{{\rm ab}}$ be 
the abelianizations of $W$ and $W'$ respectively, 
and let $\pi:W\rightarrow W^{{\rm ab}}$ 
and $\pi':W'\rightarrow {W'}^{{\rm ab}}$ be the abelianization maps.
We note that 
$W^{{\rm ab}}=(W_{\bar{S}})^{{\rm ab}}\cong\Z_2^{|\bar{S}|}$ and 
${W'}^{{\rm ab}}=(W'_{\bar{S'}})^{{\rm ab}}\cong\Z_2^{|\bar{S'}|}$.

We can obtain the following lemma 
by the same argument as the proof of \cite[Theorem~4.4]{R}.

\begin{Lemma}\label{lem3-2}
If $A$ and $B$ are subsets of $\bar{S}$ and $\pi(W_A)=\pi(W_B)$, 
then $A=B$.
\end{Lemma}

For $A\subset \bar{S}$ and $A'\subset \bar{S'}$, 
we denote $A\tau A'$ if $\pi'(\phi(W_A))=\pi'(W'_{A'})$.

We can obtain the following lemma 
by the same argument as the proof of \cite[Theorem~4.5]{R}.

\begin{Lemma}\label{lem3-3}
Let $A$ and $B$ be subsets of $\bar{S}$ and 
let $A'$ and $B'$ be subsets of $\bar{S'}$.
\begin{enumerate}
\item[(i)] If $A\tau A'$ and $B\tau A'$ then $A=B$.
\item[(ii)] If $A\tau A'$ and $A\tau B'$ then $A'=B'$.
\item[(iii)] If $A\tau A'$ and $B\tau B'$ then $(A\cap B)\tau(A'\cap B')$.
\end{enumerate}
\end{Lemma}

We obtain the following lemma from Lemmas~\ref{lem3-2} and \ref{lem3-3}.

\begin{Lemma}\label{lem3-4}
Let $A$ and $B$ be subsets of $\bar{S}$ and 
let $A'$ and $B'$ be subsets of $\bar{S'}$.
If $A\tau A'$, $B\tau B'$ and $A\subset B$, 
then $A'\subset B'$.
\end{Lemma}

\begin{proof}
Suppose that $A\tau A'$, $B\tau B'$ and $A\subset B$.
By Lemma~\ref{lem3-3}~(iii), $(A\cap B)\tau(A'\cap B')$.
Since $A\subset B$, $A\tau(A'\cap B')$.
Now $A\tau A'$.
By Lemma~\ref{lem3-3}~(ii), $A'=A'\cap B'$, 
i.e., $A'\subset B'$.
\end{proof}

A subset $T$ of $S$ is said to be {\it independent}, 
if $m(s,t)=2$ for all $s\neq t$ in $T$.
We note that 
if $T$ is an independent subset of $S$ 
then $W_T\cong \Z_2^{|T|}$.
Let $\mathcal{B}$ and $\mathcal{B}'$ be 
the sets of all maximal independent subsets of $\bar{S}$ and $\bar{S'}$, respectively.

We show the following lemma which correspond to \cite[Theorem~4.7]{R}.

\begin{Lemma}\label{lem3-5}
For each $T\in \mathcal{B}$, 
there exists a unique $T'\in \mathcal{B}'$ 
such that $T\tau T'$.
\end{Lemma}

\begin{proof}
Let $T\in \mathcal{B}$.
Then there exists $U\in\mathcal{A}$ 
such that $T\subset U\subset\bar{S}$.
By Lemma~\ref{lem2-2}, 
$\phi(W_U)=w'W'_{U'}{w'}^{-1}$ 
for some $U'\in\mathcal{A}'$ and $w'\in W'$.
Here 
$\phi:W_U\rightarrow w'W'_{U'}{w'}^{-1}$ is an isomorphism 
and $(W_U,U)$ and $(W'_{U'},U')$ are strong even.
By the proof of \cite[Theorem~4.7]{R}, 
there exists a unique independent subset $T'$ of $U'$ 
such that $T\tau T'$.
We show that $T'$ is a {\it maximal} independent subset of $\bar{S'}$.
Suppose that $T'\subset T'_0$ 
and $T'_0$ is an independent subset of $\bar{S}$.
Then by the above argument, 
there exists an independent subset $T_0$ of $\bar{S}$ 
such that $T_0\tau T'_0$.
Since $T'\subset T'_0$, $T\subset T_0$ by Lemma~\ref{lem3-4}.
Hence $T=T_0$ because $T$ is a maximal independent subset of $\bar{S}$.
By Lemma~\ref{lem3-3}~(ii), $T'=T'_0$.
Thus $T'$ is a maximal independent subset of $\bar{S'}$, 
i.e., $T'\in\mathcal{B}'$ 
which is a unique element such that $T\tau T'$.
\end{proof}

We can obtain the following lemma 
from Lemmas~\ref{lem3-3}~(iii) and \ref{lem3-5} and 
the proof of \cite[Theorem~4.8]{R}.

\begin{Lemma}\label{lem3-6}
Let $T_1,\dots,T_k\in \mathcal{A}\cup\mathcal{B}$ 
and $T'_1,\dots,T'_k\in \mathcal{A}'\cup\mathcal{B}'$ 
such that $T_i\subset \bar{S}$ and $T_i\tau T'_i$ for each $i=1,\dots,k$.
Then $|T_1\cap\dots\cap T_k|=|T'_1\cap\dots\cap T'_k|$.
\end{Lemma}

Lemma~\ref{lem3-6} implies that 
there exists a bijection $\bar{\psi}:\bar{S}\rightarrow\bar{S'}$ 
such that for each $s\in \bar{S}$ and $T\in\mathcal{A}\cup\mathcal{B}$ 
with $T\subset \bar{S}$, 
$s\in T$ if and only if $\bar{\psi}(s)\in T'$, 
where $T'$ is the element of $\mathcal{A}'\cup\mathcal{B}'$ 
such that $T\tau T'$ (cf.\ \cite{R}).
By the proof of \cite[Theorem~4.11]{R}, 
the bijection $\bar{\psi}:\bar{S}\rightarrow\bar{S'}$ 
induces an isomorphism between the Coxeter systems 
$(W_{\bar{S}},\bar{S})$ and $(W_{\bar{S'}},\bar{S'})$.

Here we note that 
we can constract $\bar{\psi}:\bar{S}\rightarrow\bar{S'}$ 
so that $\bar{\psi}(t)=t'$ 
for each $t\in \bar{S}$ and $t'\in\bar{S'}$ such that $\{t\}\tau\{t'\}$.
Indeed, suppose that $\{t\}\tau\{t'\}$ 
(such $t'$ is unique, since 
$(W_{\bar{S}},\bar{S})$ and $(W'_{\bar{S'}},\bar{S'})$ are even).
Then for $T\in\mathcal{A}\cup\mathcal{B}$ with $T\subset \bar{S}$ and 
$T'\in\mathcal{A}'\cup\mathcal{B}'$ such that $T\tau T'$, 
$t\in T$ if and only if $t'\in T'$ by Lemma~\ref{lem3-4}.

Using the above argument, we show the following.

\begin{Theorem}
The Coxeter systems 
$(W,S)$ and $(W',S')$ are isomorphic.
\end{Theorem}

\begin{proof}
We define a bijection $\psi:S\rightarrow S'$ as follows:
Let $s\in S$. 
If $s\in\bar{S}$ then we define $\psi(s)=\bar{\psi}(s)$.
Suppose that $s\in S\setminus\bar{S}$.
Then there exists a unique element 
$t\in S\setminus\{s\}$ such that $m(s,t)$ is odd.
Here we note that $m(s,u)=\infty$ for any $u\in S\setminus\{s,t\}$.
Now either $t\in\bar{S}$ or $t\not\in\bar{S}$.
We first suppose that $t\not\in\bar{S}$, 
i.e., $\{s,t\}\subset S\setminus\bar{S}$.
Then 
$\{T\in\mathcal{A}\,|\,T\cap\{s,t\}\neq\emptyset\}=\{\{s,t\}\}$.
There exists a unique $\{s',t'\}\in\mathcal{A}'$ 
such that $\phi(W_{\{s,t\}})\sim W'_{\{s',t'\}}$ by Lemma~\ref{lem2-2}.
Here $\{s',t'\}\subset S\setminus \bar{S}$ by \cite[Lemma~2.6]{H}.
We define $\psi(s)=s'$ and $\psi(t)=t'$.
Next we suppose that $t\in\bar{S}$.
Then 
$|\{T\in\mathcal{A}\,|\,T\cap\{s,t\}\neq\emptyset\}|=2$, and 
there exists a unique $T\in\mathcal{A}$ such that $t\in T\subset \bar{S}$.
By Lemma~\ref{lem2-2}, 
there exist unique $\{s',t'\}, T'\in\mathcal{A}'$ 
such that $\phi(W_{\{s,t\}})\sim W'_{\{s',t'\}}$ 
and $\phi(W_T)\sim W'_{T'}$.
The proof of \cite[Lemma~2.6]{H} implies that 
$\{s',t'\}\cap T'\neq\emptyset$ and $\phi(t)\sim s'\sim t'$.
We may suppose that $t'\in T'$.
Then $\bar{\psi}(t)=t'$ because $\{t\}\tau\{t'\}$.
We define $\psi(s)=s'$.

Then the bijection $\psi:S\rightarrow S'$ 
induces an isomorphism between the Coxeter systems 
$(W,S)$ and $(W',S')$ by the constraction of $\psi$.
\end{proof}

%

%
\end{document}